\documentclass [12pt] {amsart}

\title{Covering numbers and ``low $M^{*}$-estimate'' for quasi-convex bodies}

\author{A. E. Litvak}
\thanks{Litvak's and Milman's research was partially supported by
BSF. Research at MSRI  is supported in part by NSF grant DMS-9022140.}

\address{\hskip-\parindent A. E. Litvak,
 Department of Mathematics, Tel Aviv University, Ramat Aviv, Israel}
\email{alexandr@math.tau.ac.il}

\author{V. D. Milman}
\address{\hskip-\parindent V. D. Milman,
 Department of Mathematics, Tel Aviv University, Ramat Aviv, Israel}
\email{vitali@math.tau.ac.il}

\author{A. Pajor}
\address{\hskip-\parindent A. Pajor,
Universite de Marne-la-Vall\'ee, Equipe de Mathematiques,
                2 rue de la Butte Verte, 93166 Noisy-le-Grand Cedex, France}
\email{pajor@math.univ-mlv.fr}

\def\kkk{\qed}

\def\l{ \left}
\def\r{ \right}

\def\eps{\varepsilon}

\def\lam{\lambda}
\def\Ker{{\rm Ker}}

\def\nor{\parallel}
\def\str{\longrightarrow}
\def\tt{\theta}

\newcommand{\RR}{\mbox{\rm $~\vrule height6.5pt width0.5pt
depth0.3pt\!\!$R}}
\newcommand{\R}{\RR ^n}

\def\nx{ \nor x \nor }
\def\ny{ \nor y \nor }
\def\nxy{ \nor x+y \nor }

\def\o{ \tt }

\begin{document}


\begin{abstract}
%
  This article gives  estimates on covering numbers and
diameters of random proportional
sections and projections of  symmetric quasi-convex
bodies in $\R$.
   These results were known for the convex case and played
an essential role in development
of the theory. Because duality relations can not be applied in
the quasi-convex setting, new ingredients were
introduced that give new understanding
 for the convex case as well.
\end{abstract}

\maketitle

   \section{Introduction  and notation}

  Let $ |\cdot | $ be   
a euclidean   norm
on $\R $. Let $D$ be an ellipsoid
associated with this norm.
 Denote  $ A= \sqrt{\frac{n}{k}} \int \limits_{S^{n-1}}
\sqrt{    \sum _{i=1}^{k} x_i^2          }
\ d\sigma (x) $, where $\sigma $
is the normalized rotation invariant measure on the euclidean sphere
$ S^{n-1} $. Then $ A=A(n,k)<1$
 and $ A \str 1 $ as $ n,k \str \infty $.
 For any star-body $K$ in $\R$ define
$ M_K=  \int \limits_{S^{n-1}}
\nx d\sigma (x) $, where  $\nx$ is the gauge of $K$.
Let $M^{*}_K$ be $M_{K^0}$, where $K^0$ is the polar of $K$.
For any subsets $K_1, K_2 $ of $\R$ denote by
$N(K_1,K_2)$
 the smallest number $N$ such that there are
 N points $y_1,...,y_N$  in $K_1$ such that
$$ K_1 \subset \bigcup _{i=1}^{N} (y_i+ K_2  ) .$$

 Recall that a body $K$ is
called quasi-convex
if there is a constant $c$
such that $K+K\subset cK$, and given a $ p \in (0,1] $ a
body $K$ is called $p$-convex if for any
$\lam ,\mu >0$ satisfying
$\lam^p + \mu^p =1$ and any points $x,y\in K$
the point $\lam  x + \mu y$ belongs to $K$.
Note that for the gauge $\nor \cdot  \nor \, = \,
\nor \cdot  \nor _K  $  associated with the
the quasi-convex  ($p$-convex) body $K$ the following
inequality holds for any $x, y \in \R $
$$ \nxy \leq C \max\{ \nx  ,\ny \}
\ \ \  
 \l(  \nxy ^p \, \leq  \, 
   \nx ^p   +\ny ^p  \,   \r) . $$
 In particular every $p$-convex
body $K$  is also quasi-convex one and
$K+K\subset 2^{1/p}K$.
%
%
%
   A
 more delicate result is that for every quasi-convex body
$K$ ( $K+K\subset cK$)
there
 exists a  $q$-convex body $K_0$ such that
$K\subset  K_0 \subset 2c K$,
where $2^{1/q}=2 c$.
 This is Aoki-Rolewicz theorem
([KPR], [R], see also [K], p.47).
 In this note by a body we always mean
 a  centrally-symmetric
 compact
star-body, i.e. a body $K$ satisfying
   $tK \subset  K$ for any $t\in
  [-1, 1]$,
 however not all results use symmetry. Lemma~2,
 some analog of Lemma~4 for $p$-convex body $K$, 
and Theorem~3 hold in non-symmetric case also.




 Let us remind of the so-called
``low $M^{*}$-estimate" result. 
\newtheorem{theorem}{Theorem}
\begin{theorem}
 Let $\lam > 0 $ and $n$ be large enough.
    Let $K$ be a convex body  in $\R$ and
$\nor \cdot  \nor$ be the gauge
of $K$. Then   there exists a subspace $E$ of
 $( \R  , \nor \cdot  \nor  ) $ such that
$ \dim E = [\lam n ] $ and for any
$x \in E $ the following inequality holds
$$ \nx \geq  \frac
{ f(\lam ) }
{ 
M^{*}_K}
 | x | \ \ $$
 for some function $ f(\lam ) $, $0<\lam <1 $.

\end{theorem}
{\it Remark.} Inequality of this type  was first proved in
[M1] with very poor dependence on
$\lam $ and then improved in [M2]
to $  f(\lam )=C (1-\lam ) $. It
was later shown ([PT]), that one can take $ f(\lam )=C  \sqrt{1-\lam }$
(for different proofs see [M3] and [G]).

 By duality this theorem is equivalent
 to the following
 \medskip
 theorem.
\newline
\bf Theorem 1'    \it
Let $\lam > 0 $ and $n$ be large enough.
 For any convex body $K$ in $\R$
there exists an orthogonal projection
P of
 rank $[\lam n ]$ such that
$$ PD \subset
\frac{
 M_K}
{ f(\lam )}
 PK \ , $$
where $c$ is an absolute
   \medskip
 constant.

\rm
 In this note we will extend
 both
theorems to quasi-convex bodies. Because
duality  arguments can not be
applied to a non-convex body these two theorems
become different statements.
 Also ``"$M^*_K$" should be substituted
by an appropriate quantity not involving duality.
  Note  that by avoiding  the  use of convexity assumption
 \medskip
 we in fact simplified proof also for a convex case.
%
   \section{Main results}

 The following theorem is an extension of Theorem 1'.
\begin{theorem}
 Let $\lam > 0 $ and $n$ be large enough ($n> c/(1-\lam )^2 $).
For any $p$-convex body $K$ in $\R$
there exists an orthogonal projection
P of the rank $[\lam n ]$ such that
$$ PD \subset \frac{A_p M_K}{(1-\lam )^{1+1/p}}  PK \ , $$
where
 $ A_p = const ^{\frac{\ln (2/p)}{p}}$.

\end{theorem}

The proof of this theorem is based on the next three lemmas.
The first one was proved by
W.B.Johnson and J.Lindenstrauss
in [JL].  The second one was proved in [PT]
for convex bodies
 and is the dual form of
  Sudakov minoration theorem.
\newtheorem{lemma}{Lemma}
\begin{lemma}  There is an absolute constant $c$ such that
  if  $\eps > \sqrt{c/k} $ and $N< e^{\eps^2 k/c } $,
 then for any set of points
$y_1, ... , y_N \in \R $ and any orthogonal
projection P of  rank $k$
$$\mu \l( \{ U\in O_n \ | \ \forall j :
 A(1-\eps ) \sqrt{k/n} \    |y_j | \leq
| P U y_j | \leq A(1+\eps ) \sqrt{k/n}  \   | y_j | \} \r)
\geq $$
$$ \geq 1 - \sqrt{\pi /2} \cdot e^{-\eps ^2 k /c} $$

\end{lemma}
\begin{lemma}
 Let $K$ be a body such that $K-K \subset aK $. Then
 $N(D, tK) \leq 2 e^{2 n (a M_K / t)^2 }$.

\end{lemma}
 M.~Talagrand gave a direct simple proof
of this lemma for a convex case ([T]).
   Below, using his idea, we prove more
 general lemma for $p$-convex bodies, so we do not prove
 Lemma 2 now.
%
%
\begin{lemma}
 Let $B$ be a star-body, $K$ be a $p$-convex set,
$r \in (0,1) $, $\{x_i\} \subset rB$
 and $B \subset  \bigcup
(x_i+K  ) $.
 Then $ B\subset t_r K$, where $t_r=
\frac{1}{(1-r^p)^{1/p}} $.

\end{lemma}
{\it Proof: } Obviously $t_r=\max \{ \nx _K \  | \  x\in B \} $.
Since $B \subset  \bigcup
(x_i+K  ) $,
for any point $x$ in $B$ there are  points  $x_0$ in $rB$ and
$y$ in $K$ such that $x=x_0+y$. Then by maximality
of $t_r$ and $p$-convexity of $K$ we have
$t_r^p \leq r^p t_r^p + 1$. That proves the lemma. \kkk 
\newline
{\it Remark. } Somewhat similar argument was 
 used by N.~Kalton in dealing
 with $p$-convex sets.
 \newline \newline
%
%
%
%
{\it Proof of Theorem 2: } \newline
 Any $p$-convex body
$K$ satisfies $K-K \subset aK$ with $a=2^{1/p}$.
   By Lemma~1 and Lemma~2,
we obtain
 for $c_p \approx  2^{1/p} $
 that if
$$ c_p n \l( \frac{M_K }{t}
\r) ^{2} \leq \frac{\eps ^2 k}{c} $$
 and $\eps > \sqrt{c/k} $,  then
there exist points $x_1,...,x_N$ in $D$ and an orthogonal
projection $P$
of
 rank $k$ such that
$$ PD \subset \bigcup (Px_i + tPK) \ \mbox{ and } \ \
|Px_i| \leq (1+\eps ) \sqrt{ \frac{k}{n}  } |x_i| \ \ .$$
 Let $\lam = k/n$. Denote $r=(1+\eps ) \sqrt{ \lam } $. Lemma 3 gives us
$$
PD \subset t t_r PK \ \mbox{ for }  \
t=\frac{\sqrt{c c_p } M_K}{\eps \sqrt{\lam }}
\ \mbox{ and }
\ \eps  ^2  >  \frac{c}{\lam n} \ , \  r<1        \  \  .$$
Choose $$
\eps  = \frac{1-\sqrt{\lam }}{2 \sqrt{\lam }} \ .$$
Then for $n$ large enough
 we get
$$ PD \subset \frac{A_p M_K}{(1-\lam )^{1+1/p}}  PK \ , $$
for  $ A_p = const ^{\frac{\ln (2/p)}{p}}$ .
This completes the proof. \kkk

 Theorem 2 can be formulated in the global  \medskip  form.
\newline
  \bf
Theorem 2'   \it
 Let $K$ be a $p$-convex
body in $\R$. Then there is an orthogonal
operator $U$ such that \medskip
$$ D \subset   A^{'}_p M_K ( K+UK ) \, , \ \mbox{ where }
A^{'}_p
  = const ^{\frac{\ln (2/p)}{p}}
\  .$$

\rm
 This theorem can be proved independently, but we show
how it follows from Theorem~2.  \newline
{\it Proof of Theorem 2': }
 It follows from the proof of Theorem~2 that actually the
measure of such projections is large. So we can choose
two orthogonal subspaces $E_1, E_2$ of $\R$ such that $\dim E_1 = [n/2], \
\dim E_2 = [(n+1)/2] $ and
$$ P_i D \subset
A^{''}_p M_K
P_i K \, , $$
where $P_i$ is the projection on the space $E_i$ ($i=1,2$).
 Denote $I=id_{\R }=P_1+P_2$ and $U=P_1-P_2$. So $P_1=1/2 (I+U)$ and $P_2=1/2(I-U)$.
 Then $U$ is an orthogonal operator and for any $x\in D$ we have
$$x=P_1 x+ P_2 x \subset 1/2 \ A^{''}_p M_K (I+U) K +
1/2 \ A^{''}_p M_K (I-U) K = $$
$$
  =
 A^{''}_p M_K \frac{K+K}{2} +
  A^{''}_p M_K \frac{UK-UK}{2} = A^{'}_p M_K (K+UK)  \ .
$$
That proves Theorem~2'.   \kkk


\rm
Let us complement Lemma 2 by
 mentioning how
 covering number
$N(K, tD)$
 can be estimated.  In
convex case this estimate is given by Sudakov inequality,
using quantity
$M^{*}$. More precisely, if $K$ is
  a  convex body, then
$$N(K, tD) \leq 2 e^{c n ( M^{*}_K / t)^2 }\ .$$
 Of course,
  using  duality for a
non-convex setting leads to a weak
result, and we suggest below
a substitution for quantity $M^{*}$.



  For two quasi-convex
bodies $K, B $  define the following number
$$  M(K, B) = \frac{1}{|K|}
\int \limits_{K} \nx _B    dx \   ,$$
where $|K|$ is volume of $K$, and $\nx _B$ is the gauge of $B$.
Such numbers are considered in [MP1], [MP2] and [BMMP].

%
%
%
%
%
%
\begin{lemma}
 Let $K$ and $B$ be
 star-bodies.
 Assume  $B + B \subset aB$.  For $\o >0$ denote by
 $c_{\o }=  c_{\o } (K) $
 the best possible
 constant such that
$$
 \nor x+y \nor _K^{\o } + \nor x-y \nor _K^{\o }
  \leq 2  \cdot c_{\o }^{\o } \l( \nx _K^{\o }  +
\ny _K^{\o }  \r)
\
\mbox{ for any  } \  x, y \in \R \ .
$$
  Then
 $$N(K , tB) \leq 2 c_{\o }^n
e^{ (c n /\o  )  (a  M(K, B) / t)^{\o} }\ ,\
\mbox{ where  }  c \mbox{ is an
absolute
 constant. }$$

\end{lemma}

 Note, that in a case $\o =1$ which correspond to the general
convex case this lemma was announced in [MP2].
\newline
%
%
%
%
%
{\it Proof: }
 We follow the idea of M.Talagrand of estimating
covering numbers
in case $K=D$ ([T], see also [BLM] Proposition 4.2).
 Denote the gauge
of $K$ by $\nor \cdot \nor $ and the gauge of $B$ by $ | \cdot  | $.
  Define the measure $\mu  $ by following
$$d\mu =\frac{1}{A} e^{-\nx ^{\o }} dx \ , \mbox{ where }
%
%
\  A \   \mbox{ chosen such that } \
\int \limits_{\R } d\mu =1  \  .
$$
Let $L=
 \int \limits_{\R } |x| d\mu $. Then $\mu \{|x| \leq
2L \} \geq 1/2 $.
Let $x_1,x_2,...$ be a maximal set of points in $K$ such that
$ |x_i-x_j| \geq t$. So the sets  $x_i+\frac{t}{a} B
$ have mutually disjoint interiors.
Let $y_i =\frac{ ab}{t} x_i $ for some $b$.
Then,
by symmetry of $B$ and convexity of the function
$e^t$, we have
$$ \mu \{ y_i+bB \}
=\frac{1}{A}  \int  \limits_{bB }  e^{-\nor x+y_i \nor
 ^{\o }} dx
   = \frac{1}{2 A}
  \int  \limits_{bB }  e^{-\nor x+y_i \nor ^{\o } } +
 e^{-\nor x-y_i \nor ^{\o }}  dx
\geq  $$
 $$
 \geq  \frac{1}{A}
 \int  \limits_{bB }
e^{-\frac{1}{2}(\nor x+y_i \nor ^{\o }  +\nor x-y_i \nor ^{\o } )}  dx
%
%
%
%
%
%
%
%
%
%
%
%
  \geq \frac{1}{A}  \int  \limits_{bB }
 e^{- c_{\o }^{\o } (\nor x \nor
 ^{\o } + \nor y_i \nor
 ^{\o })} dx =  $$
$$ =
\frac{1}{A}  e^{-\nor c_{\o }y_i \nor ^{\o }}
\int  \limits_{bB }
e^{-\nor c_{\o } x \nor ^{\o }}  dx  \geq
e^{-(c_{\o } ba/t)^{\o } } c_{\o }^{-n} \mu \{ c_{\o }bB \}
$$
 Choose
$b=2L/ c_{\o }  $.
Then $ \mu \{ c_{\o }bB \} \geq 1/2 $ and, hence,
$$ N(K,tB) \leq
2    c_{\o }^n e^{(2 a L/t)^{\o } } \ . $$
 Now compute $L$.
First, the normalization constant $A$ is equal
$$A = \int \limits_{\R }
e^{-\nx ^{\o }} dx =  \int \limits_{\R }
 \int \limits_{\nx }^{\infty } (-e^{-t^{\o }})' dt dx  =
\int \limits_{0 }^{\infty } \o t^{\o -1 } e^{-t^{\o }}
\int \limits_{\nx \leq t } dx dt
= $$
$$
=
\int \limits_{\nx \leq 1 } dx
\int \limits_{0 }^{\infty } \o t^{\o +n-1} e^{-t^{\o }}dt =
|K| \cdot  \Gamma \l(1+\frac{n}{\o } \r) \ ,\
 $$
  where $ \Gamma $  is the gamma-function.
 The remaining integral is
$$
 \int \limits_{\R } |x| e^{-\nx ^{\o }} dx =
  \int \limits_{\R } |x|
 \int \limits_{\nx }^{\infty } (-e^{-t^{\o }})' dt dx
 = \int \limits_{0 }^{\infty } \o t^{\o -1 } e^{-t^{\o }}
\int \limits_{\nx \leq t } |x| dx dt
= $$
$$
=
\int \limits_{\nx \leq 1 } |x| dx
\int \limits_{0 }^{\infty } \o t^{\o +n} e^{-t^{\o }}dt =
|K|
\cdot
 M(K,B) \cdot  \Gamma \l( 1+\frac{n+1}{\o } \r) \
\ .  $$
 Using Stirling's formula we get
$$L\approx \l( \frac{n}{\o } \r) ^{1/\o } M(K, B) \ .$$
That proves the lemma. \kkk


  This lemma is an extension of Lemma 2. Indeed,
   since Euclidean space
is
a 2-smooth
 space,
 then in case $K=D$ being
an ellipsoid, we have
$c_2(D) =1$.
By direct  computation,  $M(D,B)= \frac{n}{n+1}M_B$.
 Thus,
$$N(D , tB) \leq 2  e^{ (c n   ) (  M_B
  / t)^{2} }  \, .$$

  Define the following characteristic of $K$,
 $$ \tilde M_K = \frac{1}{|K|}
\int \limits_{K} |x|   dx \,   .$$

By definition, if $K$
is a $p$-convex body, then  for any $x, y \in \R $ holds
$$
 \nor x+y \nor _K^{p } + \nor x-y \nor _K^{p }
  \leq 2  \cdot  \l( \nx _K^{p }  +
\ny _K^{p }  \r)
\,   .
$$
   So, the last lemma shows that for  $p$-convex body $K$
$$N(K , tD) \leq
2  e^{ (c n /p  )  (2
   \tilde M_K   / t)^{p} } \, .$$


 The Theorem 3 follows from this
estimate by arguments similar of that in [MP].
\begin{theorem}
 Let $\lam > 0 $ and $n$ be large enough.
 Let $K$ be a $p$-convex body  in $\R$ and
$\nor \cdot  \nor$ be the gauge
of $K$. Then   there exists  subspace $E$ of
 $( \R  , \nor \cdot  \nor  ) $ such that
$ \dim E = [\lam n ] $ and for any
$x \in E $ the following inequality holds
$$ \nx \geq  \frac{(1-\lam )^{1/2+1/p}}{a_p  \tilde M_K}
 | x | \,  ,$$
where $a_p$ depends on $p$ only (more precisely
$ a_p = const ^{\frac{\ln (2/p)}{p}}$).

\end{theorem}
{\it Proof: }
   By Lemma 4 there are points $x_1, ... , x_N$ in $K$,
such that $N<
e^{c_p n (\tilde M_K / t )^p}$ and
for any $x \in K$ there exists some $x_i$ such that $|x-x_i|<t$.
 By Lemma 1 there exists an orthogonal projection $P$
on a subspace of dimension
$\delta n$ such that
for $$c_p n \l( \frac{\tilde M_K}{t} \r) ^p <
\frac{\eps ^2
\delta n}{c} \ \mbox{ and } \
\eps > \sqrt{\frac{c}{\delta n}} $$ we have
$$b|x_i|
=(1-\eps )A \sqrt{\delta } |x_i| \leq |Px_i|
\leq (1+\eps )A \sqrt{\delta } |x_i| $$
for every $x_i$.
Let $E=\Ker P$. Then $\dim E=\lam n$, where $\lam = 1-\delta$.
 Take $x$ in $K \bigcap E $.
There is $x_i$ such that $|x-x_i|<t$.
Hence
$$|x|\leq |x-x_i| +|x_i| \leq t+\frac{|Px_i|}{b} =
t+\frac{|P(x-x_i)|}{b}
\leq
$$
$$
\leq  t+\frac{|x-x_i|}{b} \leq t(1+\frac{1}{b})
\leq    \frac{const \cdot t }{(1-\eps) \sqrt{\delta }} $$
 Therefore  for
$n$ large enough and $$t=\l(
\frac{const\cdot c_p }{\eps ^2 \delta }
\r) ^{1/p } \tilde M_K $$ we
 get
 $$ \nx \geq  \frac{ const \cdot
\eps ^2 (1-\eps ) \delta ^{1/2+1/p}}{ c_p^{1/p}  \tilde M_K}
 | x | \, .$$
To obtain our result take $\eps $, say, equal to $1/2$.
  \kkk

 As was noted in [MP2]
 in some cases
 $  \tilde M_K  << M^{*} $ and then
Theorem~3 gives better estimate
than Theorem~1 even for a convex body
(in some range of $\lam $).
 As an example,
 $K=B(l^n_1)$,
$  \tilde M_K \leq const \cdot  n^{-1/2}$,
but $ M^{*}_K \geq
const \cdot  n^{-1/2}   \medskip   (\log n )^{1/2}$.
%
%
%
%
%
%
%
%
%
%
%
%
%
%
%
%
%
%
\section{Additional remarks}

 In fact, during the proof of Theorem 2
 a more general fact was proved.
\medskip \newline
\bf Fact.
\it Let $D$ be an ellipsoid and $K$ be a $p$-convex body. Let
 $$N(D,  K)\leq e^{\alpha  n} \ . $$   Denote for an integer $1\leq k \leq n$
the ratio $\lam =k/n $.
Then for some absolute constant $c$ and
$$\gamma =c \sqrt{\alpha }  \ ,\  \ k\in (\gamma ^2 n, (1-2 \gamma )^2 n) $$
there exists an orthogonal projection $P$ of rank $k$ such that
 $$
    \l( p(1- \sqrt{\lam })/2 \r) ^{1/p}
PD \subset PK  \ .
                                                               $$
\rm

 In terms of  entropy numbers this means
$$
   \frac{\l( p(1-\sqrt{k/n} )/2 \r) ^{1/p}}{e_k (D,K) } PD \subset PK  \ ,
                                                                           $$
where $e_k (D, K)= \inf \{ \eps >0 \  | \  N(D, \eps K ) \leq 2^{k-1} \} $ .

  It is worth to point out that Theorem~2 can be obtained from this
results.

  We thank E.~Gluskin for his remarks on the first
draft of this note.

\section*{References}
{\footnotesize

\begin{description}
\smallbreak
\item[{[BLM]}]
 Bourgain, J. ; Lindenstrauss, J. ; Milman, V. {\em
 Approximation of zonoids by zonotopes.}
Acta Math. 162 (1989), no. 1-2, 73--141.
\smallbreak
\item[{[BMMP]}]
  Bourgain, J. ; Meyer, M. ; Milman, V. ; Pajor, A. {\em
On a geometric inequality.}
Geometric aspects of functional analysis (1986/87), 271--282,
 Lecture Notes in Math., 1317,
Springer, Berlin-New York, 1988.
\smallbreak
\item[{[JL]}]
 Johnson, W. B.; Lindenstrauss, J.
{\em  Extensions of Lipschitz mappings into a
Hilbert space.}  Conference in modern analysis and probability
(New Haven, Conn., 1982), 189--206,
\smallbreak
\item[{[G]}]
Gordon, Y.
{\em On Milman's inequality and random subspaces which escape through a
mesh in $\R $.}
Geometric aspects of functional analysis (1986/87), 84--106, Lecture Notes
in Math., 1317, Springer, Berlin-New York, 1988.{\hfuzz=7pt\par}
\smallbreak
\item[{[KPR]}]    N.J.Kalton,  N.T.Peck, J.W.Roberts,
 {\em  An $F$-space sampler},  London
Mathematical Society Lecture Note Series, 89,
Cambridge University Press, Cambridge and New York,
1984.
\smallbreak
\item[{[K]}] K\"onig, H. {\em Eigenvalue Distribution of Compact Operators },
Birkh\"auser, 1986.
\smallbreak
\item[{[M1]}]  Milman, V. D. Almost Euclidean quotient spaces of subspaces of a
finite-dimensional normed space.
Proc. Amer. Math. Soc. 94 (1985), no. 3, 445--449.
\smallbreak
\item[{[M2]}]
 Milman, V. D. {\em
Random subspaces of proportional dimension of finite-dimen\-sional
normed spaces: approach through
the isoperimetric inequality.} Banach spaces (Columbia, Mo., 1984),
106--115,
Lecture Notes in Math., 1166, Springer, Berlin-New York, 1985.
\smallbreak
\item[{[M3]}]
 Milman, V. D. {\em
A note on a low $M\sp *$-estimate.}  Geometry of Banach spaces (Strobl,
1989), 219--229, London Math. Soc. Lecture Note Ser.,
158, Cambridge Univ. Press, Cambridge,
1990.
\smallbreak
\item[{[MP1]}]
 Milman, V. D. ; Pajor, A. {\em
Isotropic position and inertia ellipsoids and zonoids of the
unit ball of a normed $n$-dimensional space.}
Geometric aspects of functional analysis (1987--88),
64--104, Lecture Notes in Math., 1376,
Springer, Berlin-New York, 1989.
\smallbreak
\item[{[MP2]}] Milman, V. ; Pajor, A. {\em
Cas limites dans des in\'egalit\'es du type de Khinchine et
applications g\'eom\'etriques.} (French)
[Limit cases of Khinchin-type inequalities and some geometric
applications] C. R. Acad. Sci.
Paris S\'er. I Math. 308 (1989), no. 4, 91--96.
\smallbreak
\item[{[PT]}]
 Pajor, A. ; Tomczak-Jaegermann, N. {\em
Subspaces of small codimension of
finite-dimensional Banach spaces. }
Proc. Amer. Math. Soc. 97 (1986), no. 4, 637--642.
\smallbreak
\item[{[R]}]
 Rolewicz, S. {\em Metric linear spaces.}  Monografie Matematyczne, Tom. 56.
[Mathematical Monographs, Vol. 56] PWN-Polish Scientific Publishers,
Warsaw, 1972.

\end{description}
}

\end{document}